\documentclass[11pt]{article}
\usepackage{amsmath, amssymb, amsfonts, amsthm, graphicx}
\theoremstyle{plain}

\theoremstyle{definition}

\begin{document}
\title{A method for recursively generating sequential rational approximations to $\sqrt[n]{k}$}
\author{Joe Nance\\
  \texttt{nance2uiuc@gmail.com}\\
  Department of Mathematics \\
  University of Illinois at Urbana-Champaign\\}
\date{\today}
\maketitle

\begin{abstract}
The goal of this paper is to derive a simple recursion that generates a sequence of fractions approximating $\sqrt[n]{k}$ with increasing accuracy. The recursion is defined in terms of a series of first-order non-linear difference equations and then analyzed as a discrete dynamical system. Convergence behavior is then discussed in the language of initial trajectories and eigenvectors, effectively proving convergence without notions from standard analysis of infinitesimals.
\end{abstract}
\clearpage
\section{Introduction and motivation}
Consider for a moment the simple recursion
\begin{equation}\label{eq:1}
\frac{x_{t+1}}{y_{t+1}}=\frac{x_t + 2 y_t}{x_t + y_t}.
\end{equation}
Choose initial values to set it marching on its way towards $\sqrt{2}$. If we choose initial values $x_0=y_0=1$, the recursion (\ref{eq:1}) gives a sequence of fractions approximating $\sqrt{2}$ whose behavior is summarized in the table below.
\begin{center}
\begin{tabular}{| l| c| l |}
\hline
t & $x_t/y_t$ & $\approx$ \\
\hline
0 & $1/1$ & 1 \\
1 & $3/2$ & 1.5 \\
2 & $7/5$ & 1.4 \\
3 & $17/12$ & 1.4167... \\
4 & $41/29$ & 1.41379... \\
5 & $99/70$ & 1.41429... \\
\hline
\end{tabular}
\end{center}

\noindent For reassurance that the recursion generates a sequence that does infact converge to $\sqrt{2}$, employ the following analysis: For some sequence $a_t$,\newline
if $a_t\rightarrow L$, then $a_{t+1}\rightarrow L$ as well. So suppose that the recursion (\ref{eq:1}) has limit $L$.  Multiply the top and bottom of the right hand side of (\ref{eq:1}) by $1/y_t$ and we get
\[
\frac{x_{t+1}}{y_{t+1}}=\frac{x_t + 2 y_t}{x_t + y_t}\ \frac{\frac{1}{y_t}}{\frac{1}{y_t}}.
\]
 Using the previous fact about limits, we have $$L=\frac{L+2}{L+1}.$$ This gives $L=\pm \sqrt{2}$. In this analysis, we have picked up the unsettling possibility of this recursion generating a sequence of fractions converging to $-\sqrt{2}$. Discussion of which initial values generate such a sequence is withheld momentarily.

\section{$\sqrt{k}$ Recursion}
It is convenient to consider recursion (\ref{eq:1}) as an action on a system of first order linear difference equations given by
\begin{equation}\label{eq:2}
\begin{cases}
x_{t+1} = x_t+2y_t, \\
y_{t+1} = x_t + y_t,\end{cases}
\end{equation}
which is a discrete dynamical system, or ``DDS". Clearly, if we replace $2$ by a positive integer $k\in \mathbb{Z}^+$ we obtain a recursion similar in structure to recursion (\ref{eq:1}) converging to $\pm \sqrt{k}$ for any initial values $x_0,y_0$. This recursion is
\begin{equation}\label{eq:6}
\frac{x_{t+1}}{y_{t+1}}=\frac{x_t + k y_t}{x_t + y_t}.
\end{equation}
This gives a corresponding DDS
\begin{equation}\label{eq:3}
\begin{cases}
x_{t+1} = x_t+k y_t, \\
y_{t+1} = x_t + y_t,\end{cases}
\end{equation}
which can be represented in matrix form as
\begin{equation}\label{eq:4}
\begin{pmatrix}
x_{t+1} \\
y_{t+1} \end{pmatrix}
=
\begin{pmatrix}
1 & k \\
1 & 1 \end{pmatrix}
\begin{pmatrix}
x_t \\
y_t \end{pmatrix}.
\end{equation}
Any term in a sequence generated by (\ref{eq:4}) is generalized as
\begin{equation}\label{eq:5}
\begin{pmatrix}
x_t \\
y_t\end{pmatrix}
={\begin{pmatrix}
1 & k \\
1 & 1\end{pmatrix}}^t
\begin{pmatrix}
x_0 \\
y_0\end{pmatrix}.
\end{equation}
Recursion (\ref{eq:6}) is recovered by taking ratios of terms with equal indices.

\section{$\sqrt[n]{k}$ Recursion}
A natural question at this point would be ``Does there exist a structurally simple recursion similar to (\ref{eq:6}) generating a sequence of fractions approximating $\sqrt[n]{k}$ ?", to which the answer is ``kind of". In the spirit of the previous analysis, start out with \[L=\sqrt[n]{k}\] so that \[L^n=k.\] Add $L$ to both sides to get \[L^n + L = L+k,\] then factor out an $L$: \[L(L^{n-1}+1)=L+k.\] Now divide both sides by $L^{n-1}+1$ to get \[L=\frac{L+k}{L^{n-1}+1}.\] Use the fact about sequences to obtain
\begin{equation}\label{eq:7}
\frac{x_{t+1}}{y_{t+1}}=\frac{x_t + k y_t}{\frac{x^{n-1}_t}{y^{n-2}_t} + y_t},
\end{equation} which has the corresponding DDS
\begin{equation}\label{eq:8}
\begin{cases}
x_{t+1} &=  x_t+k y_t \\
y_{t+1} &=  \frac{x_t^{n-1}}{y_t^{n-2}} + y_t.\end{cases}
\end{equation} This latter system is not representable as a simple $2\times 2$ matrix with real entries for arbitrary $n$. This is because there are now \emph{three} different terms in this system: $x_t$, $y_t$, and $x_t^{n-1}/y_t^{n-2}$ which is non-linear.

\section{Convergence behavior of the DDS given by $\Pi_n$}
We are now in a position to analyze (\ref{eq:8}) using matrix methods. Suppose there \emph{did exist} an $n\times n$ matrix, call it $\Pi_n$, such that (\ref{eq:8}) could be represented as
\begin{equation}\label{eq:11}
\begin{pmatrix}
\vec{x}_{t,1} \\
\vec{x}_{t,2} \\
\vdots \\
\vec{x}_{t,n}\end{pmatrix}=\Pi_n^t
\begin{pmatrix}
\vec{x}_{0,1} \\
\vec{x}_{0,2} \\
\vdots \\
\vec{x}_{0,n}\end{pmatrix},
\end{equation} where $\vec{x}_{a,b}$ denotes the $b$th entry of the vector corresponding to $a$.
Or equivalently, $\vec{R}_t=\Pi_n^t \vec{R}_0$. But what would the matrix $\Pi_n$ look like? Consider the properties of (\ref{eq:8}) which $\Pi_n$ must capture in order to faithfully represent convergence behavior of (\ref{eq:7}). First, recognize that the action on (\ref{eq:8}) needed to arrive at (\ref{eq:7}) is taking ratios of terms with equal indices, which is equivalent to taking the ratio of successive entries of a vector $\vec{R}_t$. The entries in an arbitrary vector $\vec{R}_t$ must tend toward those of the dominant eigenvector of $\Pi_n$, $\vec{\lambda}_d$. So the ratio of any two successive entries in $\vec{\lambda}_d$ should be equal to $\sqrt[n]{k}$. That is, $\vec{\lambda}_{d,i}/\vec{\lambda}_{d,i+1}=\sqrt[n]{k}$ where $1\leq i\leq n-1$ and $\vec{\lambda}_{d,i}$ is the $i$th entry of the dominant eigenvector of $\Pi_n$. Such an eigenvector looks like
\[
\vec{\lambda}_d=\begin{pmatrix}
\sqrt[n]{k}^{n-1} \\
\sqrt[n]{k}^{n-2} \\
\vdots \\
 1\end{pmatrix}.
 \]
Recovery of the dominant eigenvalue from the dominant eigenvector is had by solving for $\lambda_d$ in $(\Pi_n-\lambda_d \textit{I}_n)\vec{\lambda}_d=\vec{0}$ where $\vec{0}$ is the column vector consisting entirely of zeros. This calculation yields $\lambda_d=1+\sqrt[n]{k}$.

Now suppose $\Pi_n$ is diagonalizable. Then $\Pi_n$ admits a basis for $\mathbb{R}^n$ consisting entirely of eigenvectors of $\Pi_n$. So any initial vector $\vec{R}_0$ in $\mathbb{R}^n$ can be written as a linear combination of eigenvectors $\vec{R}_0=\sum\limits_{i=1}^n c_i\vec{\lambda}_i$.

Applying $\Pi_n$ to our initial vector $\vec{R}_0$,
\begin{equation*}\begin{split}
\vec{R}_1 = \Pi_n \vec{R}_0 &= \Pi_n\displaystyle\sum\limits_{i=1}^n c_i\vec{\lambda}_i \\
&= c_1\Pi_n\vec{\lambda}_1+\cdots+c_n\Pi_n\vec{\lambda}_n \\
&= c_1\lambda_1\vec{\lambda}_1+\cdots+c_n \lambda_n \vec{\lambda}_n.\end{split}
\end{equation*}
Applying $\Pi_n$ again,
\begin{equation*}\begin{split}
\vec{R}_2=\Pi_n^2\vec{R}_0 = \Pi_n(\Pi_n\vec{R}_0) &=\Pi_n(c_1\lambda_1\vec{\lambda}_1+\cdots+c_n\lambda_n\vec{\lambda}_n) \\
&=c_1\lambda_1\Pi_n\vec{\lambda}_1+\cdots+c_n\lambda_n\Pi_n\vec{\lambda}_n \\
&=c_1\lambda_1^2{\vec{\lambda}_1}+\cdots+c_n\lambda_n^2{\vec{\lambda}_n}.\end{split}\end{equation*}
We can see the pattern now,
\begin{equation}\label{eq:9}
\vec{R}_t=\Pi_n^t\vec{R}_0 = c_1{\lambda_1}^t\vec{\lambda}_1+\cdots+c_n{\lambda_n}^t\vec{\lambda}_n=\displaystyle\sum\limits_{i=1}^n c_i\lambda_i^t\vec{\lambda}_i.
\end{equation}
Verify that indeed,
\begin{equation*}
\vec{R}_t=\Pi_n^t\vec{R}_0 =\begin{pmatrix}
c_1{\lambda_1}^t\vec{\lambda}_{1, 1}+\cdots+c_n{\lambda_n}^t\vec{\lambda}_{n, 1} \\
c_1{\lambda_1}^t\vec{\lambda}_{1, 2}+\cdots+c_n{\lambda_n}^t\vec{\lambda}_{n, 2} \\
\vdots \\
c_1{\lambda_1}^t\vec{\lambda}_{1, n}+\cdots+c_n{\lambda_n}^t\vec{\lambda}_{n, n}\end{pmatrix},
\end{equation*}
where $\vec{\lambda}_{a,b}$ denotes the $b$th entry of the eigenvector corresponding to $\lambda_a$.

Note that with increasing $t$, only one of the $c_i\lambda_i^t\vec{\lambda}_{i,l}$ terms becomes the dominant term. The dominant term is the one involving the dominant eigenvalue, $\lambda_d$ and the entries of its corresponding eigenvector, $\vec{\lambda}_{d,i}$. Since contributions of the other terms become negligible in the limiting quotient, we can make the following statement:
\begin{equation*}
\vec{R}_t=\Pi_n^t \vec{R}_0\approx c\lambda_d^t\vec{\lambda}_d=c(1+\sqrt[n]{k})^t\begin{pmatrix}
\left(\sqrt[n]{k}\right)^{n-1} \\
\left(\sqrt[n]{k}\right)^{n-2} \\
\vdots \\
1\end{pmatrix}\mbox{ for } t\gg 1.
\end{equation*}
It follows that
\begin{equation}\label{eq:18}
\begin{split}
\frac{\vec{R}_{t,i}}{\vec{R}_{t,i+1}}&=
\frac{c(1+\sqrt[n]{k})^t\vec{\lambda}_{d,i}+\cdots+c_n{\lambda_n}^t\lambda_{n,i}}
{c(1+\sqrt[n]{k})^t\vec{\lambda}_{d,i+1}+\cdots+c_n{\lambda_n}^t\lambda_{n,i+1}} \\
&\approx\frac{c(1+\sqrt[n]{k})^t\vec{\lambda}_{d,i}}{c(1+\sqrt[n]{k})^t\vec{\lambda}_{d,i+1}}=\frac{\vec{\lambda}_{d,i}}{\vec{\lambda}_{d,i+1}}=\sqrt[n]{k},
\end{split}
\end{equation}
for $(1\leq i\leq n-1)$ and $t\gg1$.

When $n=2$ (finding square roots), let $c=0$ and notice that (\ref{eq:18}) is approximately $-\sqrt{k}$. To satisfy our curiosity from Section 1, if we are to have a sequence generated by (\ref{eq:2}) converging to $-\sqrt{k}$ then the appropriate initial values are ones which are components of some multiple of the second eigenvector, not a linear combination of the dominant eigenvector and the other. This way, ratios of successive entries of an evolving vector equal the slope of the second eigenvector, $-\sqrt{k}$. If we restrict our choices of initial values to $\mathbb{Q}^2$, then we do not run into this problem of multiple limits.

\section{Derivation of the matrix $\Pi_n$}
In the preceding section we showed that if (\ref{eq:8}) can be represented as (\ref{eq:11}), then  $\lim_{t\to+\infty} x_{t,i}/x_{t,i+1}=\sqrt[n]{k}$, as desired. It only remains to find the exact form for $\Pi_n$. In the analysis above, we showed that long-term time evolution of an initial vector depends heavily on ``hitting" $\vec{\lambda}_d$ with its corresponding eigenvalue $\lambda_d$, which in turn depends on left-multiplying $\vec{R}_0$ by $\Pi_n$. Carrying this calculation out gives
\[
\lambda_d \vec{\lambda}_d=\left(1+\sqrt[n]{k}\right)\begin{pmatrix}
{\left(\sqrt[n]{k}\right)}^{n-1} \\
{\left(\sqrt[n]{k}\right)}^{n-2} \\
{\left(\sqrt[n]{k}\right)}^{n-3} \\
{\left(\sqrt[n]{k}\right)}^{n-4} \\
\vdots \\
1\end{pmatrix}
=
\begin{pmatrix}
{\left(\sqrt[n]{k}\right)}^{n-1} \\
{\left(\sqrt[n]{k}\right)}^{n-2} \\
{\left(\sqrt[n]{k}\right)}^{n-3} \\
{\left(\sqrt[n]{k}\right)}^{n-4} \\
\vdots \\
1\end{pmatrix}
+
\begin{pmatrix}
k \\
{\left(\sqrt[n]{k}\right)}^{n-1} \\
{\left(\sqrt[n]{k}\right)}^{n-2} \\
{\left(\sqrt[n]{k}\right)}^{n-3} \\
\vdots \\
\sqrt[n]{k}\end{pmatrix}.
\]

Evidently, $\Pi_n$ is such that when an initial vector $\vec{R}_0$ is left multiplied by it, returned is the dominant eigenvector plus another vector, $\vec{v}$, along with other negligible terms in the limit that $t\gg 1$. This means that $\Pi_n$ must be the sum of two matrices acting on the linear combination of eigenvectors that comprises $\vec{R}_0$. \[\Pi_n \vec{R}_0 = \left(\textit{I}_n+\pi_n\right)\left(\lambda_d \vec{\lambda}_d+\cdots+\lambda_n \vec{\lambda}_n\right).\] where $\textit{I}_n$ is the $n$-dimensional identity matrix. But what is $\pi_n$? Consider what action is taken by $\pi_n$ to return $\vec{v}$ from $\vec{\lambda}_d$. Apparently, $\pi_n$ is an $n\times n$ matrix such that when it left multiplies a column vector, it has the effect of permuting entries by one place in a cyclic manner while scaling by a factor of $\sqrt[n]{k}$. By inspection, we see that
\[
\pi_n=\begin{pmatrix}
0 & 0 & 0 & \cdots  & 0 & k \\
1 & 0 & 0 & \cdots  & 0 & 0 \\
0 & 1 & 0 & \cdots  & 0 & 0 \\
0 & 0 & 1 & \cdots  & 0 & 0 \\
\vdots & \vdots  & \vdots  & \ddots  & \vdots  & \vdots  \\
0 & 0 & 0 & \cdots  & 1 & 0\end{pmatrix}.\] The exact form of $\Pi_n$ is given by
\[
\begin{split}
\Pi_n&=\textit{I}_n+\pi_n \\
&=\begin{pmatrix}
1 & 0 & 0 & \cdots & 0 & 0 \\
0 & 1 & 0 & \cdots & 0 & 0 \\
0 & 0 & 1 & \cdots & 0 & 0 \\
0 & 0 & 0 & \cdots & 1 & 0 \\
\vdots & \vdots  & \vdots  & \ddots  & \vdots  & \vdots  \\
0 & 0 & 0 & \cdots & 0  & 1\end{pmatrix}
+\begin{pmatrix}
0 & 0 & 0 & \cdots  & 0 & k \\
1 & 0 & 0 & \cdots  & 0 & 0 \\
0 & 1 & 0 & \cdots  & 0 & 0 \\
0 & 0 & 1 & \cdots  & 0 & 0 \\
\vdots & \vdots  & \vdots  & \ddots  & \vdots  & \vdots  \\
0 & 0 & 0 & \cdots  & 1 & 0\end{pmatrix} \\
&=\begin{pmatrix}
1 & 0 & 0 & \cdots  & 0 & k \\
1 & 1 & 0 & \cdots  & 0 & 0 \\
0 & 1 & 1 & \cdots  & 0 & 0 \\
0 & 0 & 1 & \cdots  & 0 & 0 \\
\vdots & \vdots  & \vdots  & \ddots  & \vdots  & \vdots  \\
0 & 0 & 0 & \cdots  & 1 & 1\end{pmatrix}.\end{split}
\]

\section{Algebraic properties of the matrix $\Pi_n$}
We must still confirm diagonalizability of $\Pi_n$ since most of our case depends upon this property of $\Pi_n$. The characteristic polynomial, $P(\lambda)$, of $\Pi_n$ can be found by computing $\det(\Pi_n-\lambda \textit{I}_n)$ by expanding in minors along the top row, giving $P(\lambda)=(1-\lambda)^n+(-1)^{n+1}k$. The eigenvalues are had by solving $P(\lambda)=0$ giving $\lambda_j=1-\sqrt[n]{k}\ e^{\textit{i}J\pi/n}$ where $1\leq j \leq n$ and
$J=2j\mbox{ when } n \text{ is even and } 2j+1\text{ when }n\text{ is odd.}$ The largest of these eigenvalues is $\lambda_d=\lambda_{\frac{n}{2}}=\lambda_{\frac{n-1}{2}}=1+\sqrt[n]{k}$ as desired. These $n$ distinct eigenvalues give $n$ distinct eigenvectors given by
\[
\vec{\lambda}_j=\begin{pmatrix}
\frac{k}{\lambda_j-1} \\
 \frac{k}{(\lambda_j-1)^2} \\
  \frac{k}{(\lambda_j-1)^3} \\
   \frac{k}{(\lambda_j-1)^4} \\
    \vdots \\
    1\end{pmatrix},\]
     the largest of which corresponds to the eigenvalue $\lambda_d$; this vector is
\[
\vec{\lambda}_d=\begin{pmatrix}
\left(\sqrt[n]{k}\right)^{n-1} \\
 \left(\sqrt[n]{k}\right)^{n-2} \\
  \left(\sqrt[n]{k}\right)^{n-3} \\
   \left(\sqrt[n]{k}\right)^{n-4} \\
    \vdots \\
     1\end{pmatrix}.
     \] Since $\Pi$ gives $n$ distinct eigenvalues and eigenvectors, $\Pi_n$ is diagonalizable. We are now in the position to make the generalization of (\ref{eq:7}) as follows
\begin{equation}\label{eq:12}
\Pi_n^t \vec{R}_0 =\vec{R}_t \text{ gives rise to } \displaystyle\lim_{t\to+\infty}\frac{\vec{R}_{t,i}}{\vec{R}_{t,i+1}}=\sqrt[n]{k}
\end{equation}
for $1\leq i \leq n-1$ and $t\in \mathbb{Z}$.
Equivalently,
\begin{equation*}
\begin{pmatrix}
1 & 0 & 0 & \cdots  & 0 & k \\
1 & 1 & 0 & \cdots  & 0 & 0 \\
0 & 1 & 1 & \cdots  & 0 & 0 \\
0 & 0 & 1 & \cdots  & 0 & 0 \\
\vdots & \vdots  & \vdots  & \ddots  & \vdots  & \vdots  \\
0 & 0 & 0 & \cdots  & 1 & 1\end{pmatrix}^t
\begin{pmatrix}
x_1 \\
x_2 \\
x_3 \\
x_4 \\
\vdots \\
x_n\end{pmatrix}
=
\begin{pmatrix}
x^\prime_1 \\
x^\prime_2 \\
x^\prime_3 \\
x^\prime_4 \\
\vdots \\
x^\prime_n\end{pmatrix}
\end{equation*}
so that $\displaystyle\lim_{t\to+\infty}\frac{x^\prime_i}{x^\prime_{i+1}}=\sqrt[n]{k}$.

This result effectively fulfils the goal of the paper which was to derive a simple recursion that generates a sequence of fractions approximating $\sqrt[n]{k}$ with increasing accuracy.

\section{Computation}
The result of the previous section satisfies the technical goal of this paper, but it is left to the reader to judge the practicality of this result. Accuracy of an approximation depends on taking powers of an $n\times n$ matrix. This tedious task can be tiresome for even relatively small powers of $n$ and $t$. So where do we look to find aid in this computation? One could certainly start with the Cayley-Hamilton Theorem, which states that every $n\times n$ matrix over a commutative ring satisfies its own characteristic equation, $P(\lambda)=$ $\det(\lambda \textit{I}_n-A).$ Applying this theorem to $\Pi_n$ gives \[P(\Pi_n)=(I_n-\Pi_n)^n + (-1)^{n+1} k \ I_n=\mathbf{0},\] where $\textbf{0}$ is the $n\times n$ matrix consisting entirely of zeros. The binomial theorem then gives \[\sum_{i=0}^n \frac{n!(-1)^i}{i!(n-i)!}\Pi_n^i I_n = (-1)^n k\ I_n.\] Solving then for $\Pi_n^n$ gives
\begin{equation}\label{eq:14}
\Pi_n^n=\sum_{i=1}^{n-1} \frac{n!(-1)^{n-1-i}}{i!(n-i)!}\Pi_n^i +\left[(-1)^{n-1}+k\right]I_n.
\end{equation}
This is an explicit equation expressing $\Pi_n^n$ in terms of lower powers of $\Pi_n$ and $\textit{I}_n$. It is useful because if one is able to calculate powers of $\Pi_n$ up to and including $\Pi_n^{n-1}$, then one is able to generate arbitrarily large powers of $\Pi_n$ iteratively which then can be used to generate arbitrarily close approximations  to $\sqrt[n]{k}$.
\newline

Let's take a look at the $n=2$ case. Equation (\ref{eq:14}) gives
\begin{equation}\label{eq:15}
\Pi_2^2=2\Pi_n+(k-1)\textit{I}_2,
\end{equation}
which is an explicit expression of $\Pi_2$ in first powers of $\Pi_2$ and $\textit{I}_2$. Because no higher powers of $\Pi_2$ need to be calculated to arrive at (\ref{eq:15}), arbitrary integer powers of $\Pi_2$ are gotten with ease from iterative multiplication and substitution of powers of $\Pi_2$. This gives rise to a Fibonacci-like sequence in the exponents of $\Pi_2$:

\begin{equation*}\begin{split}
\Pi_2^2 &=2\Pi_n+(k-1)\textit{I}_2\\
\Pi_2^3  &= \Pi_2 \Pi_2^2 \\
&= 2\Pi_2^2 + (k-1)\Pi_2\\
&= (k+3)\Pi_2+2(k-1)\textit{I}_2\\
\Pi_2^5  &= \Pi_2^3 \Pi_2^2 \\
&= (k^2+10k+5)\Pi_2+4(k^2-1)\textit{I}_2\\
&\vdots  \\
\Pi_2^{F_i} &= \Pi_2^{F_{i-1}} \Pi_2^{F_{i-2}}.\end{split}
\end{equation*}
 The reader is encouraged to try this for the $n=3,4,5,...$ cases to see that once harrowing computations are done to make $\Pi_n^{n-1}$  known, precise approximate computation soon follows.

\end{document}